\documentclass[11pt]{article}
\usepackage[font=small,format=plain,labelfont=bf,up]{caption}

\usepackage[a4paper,nohead,ignorefoot,%
twoside,scale=0.875,bindingoffset=1.5cm]{geometry}

\usepackage{amssymb,amsfonts,amsmath,amsthm}

\usepackage[]{graphicx}
\usepackage{color}

\graphicspath{{}}

\definecolor{colorGreen}{rgb}{0.,0.67,0}
\definecolor{colorRed}{rgb}{0.9,0,0}
\definecolor{colorBlue}{rgb}{0.,0.,0.67}

\newcommand{\BMHC}{}
\newcommand{\EMHC}{}

\usepackage{float}
\usepackage{mathrsfs}
\usepackage{marginnote}

\newtheorem*{Cor*}{Corollary}

\newtheorem*{Thm*}{Theorem}
\newtheorem*{Observation*}{Observation}

\newtheorem*{Prop*}{Proposition}
\theoremstyle{definition}

\newtheorem*{Rmk*}{Remark}

\newtheorem*{Ex*}{Example}

\newtheorem*{Qu*}{Question}

\newcommand{\R}{\mathbb{R}}
\newcommand{\om}{\omega}

\newcommand{\beq}{\begin{equation}}
\newcommand{\beqn}{\begin{equation}\nonumber}
\newcommand{\eeq}{\end{equation}}

\newcommand{\bea}{\begin{equation}\begin{aligned}}
\newcommand{\bean}{\begin{equation}\begin{aligned}\nonumber}
\newcommand{\eea}{\end{aligned}\end{equation}}

\numberwithin{equation}{section}

\definecolor{Urs}{rgb}{0,.7,0}
\definecolor{Peter}{rgb}{0,0,1}
\definecolor{red}{rgb}{1,0,0}

\renewcommand\footnotemark{}
\begin{document}
%
%
% -----------------------------------------------------------------------------
% - Title information
% -----------------------------------------------------------------------------
%
%
\title{Numerical simulations of magnetic billiards in a convex domain in $\mathbb{R}^2$}
\date{\today}
\author{%
Peter Albers\footnote{ %
Institutes for Applied and Pure Mathematics, University of M{\"u}nster, Germany
} 
\and
Gautam Dilip Banhatti\thanks{\emph{E-mail: }%
{\tt peter.albers/g\_banh02/michael.herrmann@uni-muenster.de}
} %
\and 
Michael Herrmann
%\thanks{%
%Numerical data available at {\tt www....}
%} %
%
 %
} %
\maketitle
%
%
% -----------------------------------------------------------------------------
% - Abstract
% -----------------------------------------------------------------------------
%
%
%
\begin{abstract}
We present numerical simulations of magnetic billiards inside a convex domain in the plane.
\end{abstract}
%%
%%
% -----------------------------------------------------------------------------
% - MSC and keywords
% -----------------------------------------------------------------------------
%
%%\quad\newline\noindent%
%%\begin{minipage}[t]{0.15\textwidth}%
%%   Keywords: 
%%\end{minipage}%
%%\begin{minipage}[t]{0.8\textwidth}% 
%%\end{minipage}%
%%% 
%%\medskip
%%\newline\noindent
%%\begin{minipage}[t]{0.15\textwidth}%
%%   MSC (2010): %
%%\end{minipage}%
%%\begin{minipage}[t]{0.8\textwidth}%
%%\end{minipage}%

% -----------------------------------------------------------------------------
% - Table of contents
% -----------------------------------------------------------------------------
%
%\setcounter{tocdepth}{3}
%\setcounter{secnumdepth}{3}{\scriptsize{\tableofcontents}}
%
%
%
% -------------------------------------------------------------------------------------
\section{Introduction}
% -------------------------------------------------------------------------------------
%

In this article we present some numerical simulations of magnetic billiards inside a convex domain in $\R^2$. \BMHC Classical \EMHC Billiards is a simple dynamical system which shows up in \BMHC various branches of mathematics. \EMHC It is extremely well studied with deep results and many open questions at the same time.

Robnik and Berry, \cite{Robnik_Berry_Classical_billiards_in_magnetic_fields}, were the first to study numerically magnetic billiards in the plane. It seems that this article is still a main source for numerical results of classical magnetic billiards, in particular in ellipses, see for instance page 1 in \cite{Bialy_Mironov_Algebraic_non_integrability_of_magnetic_billiards}. Newer sources for numerical results are for instance \cite{Meplan_Brut_Gignoux_Tangent_map_for_classical_billiards_in_magnetic_fields,Berglund_Kunz_Integrability_and_ergodicity_of_classical_billiards_in_a_magentic_field}.
The quantum mechanical analogue seems to be much more studied but doesn't concern us here.
 
The purpose of this article is to provide a more detailed numerical study of classical magnetic billiards inside a convex domain (mostly ellipses). We do not claim any originality and consider this purely as a service to the community. We used Matlab for the numerical computations
and Mathematica to illustrate the results.

\section{Magnetic billiards inside a convex domain}

We briefly describe the dynamical system. For that let $\Sigma\subset \R^2$ be a curve bounding a strictly convex domain $T$. For our purposes being a curve means that $\Sigma$ is a smooth embedded compact 1-manifold without boundary. This assumptions are very restrictive and, for instance, exclude billiards in polygons. A much more general account and set-up can be found in the book by Tabachnikov, \cite{Tabachnikov_Geometry_and_billiards}. The domain $T$ bounded by $\Sigma$ is called the table. 

\begin{figure}[ht!]
\def\svgwidth{80ex}
\centering{ %
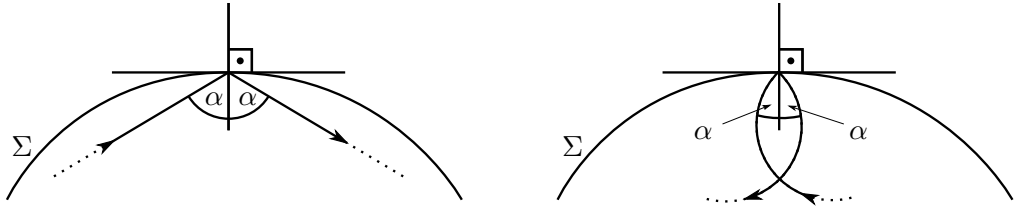
} %
\caption{Billiard reflection without (left) and with (right) magnetic field}\label{fig:billiard}
\end{figure}

Non-magnetic billiard inside $T$ describes the motion of a free particle which undergoes elastic reflection at the boundary, that is, the point particle moves with constant speed on a straight line until it hits the boundary. Then this straight line is reflected according to the law 'angle of incidence = angle of reflection', \BMHC i.e., \EMHC the tangential component of the velocity is kept whereas the normal component is flipped, see Figure \ref{fig:billiard}. 

\BMHC For magnetic billiards, a charged particle moves \EMHC in a constant magnetic field which is perpendicular to the plane containing the table. Thus the particle moves on a circle instead of a straight line. The reflection law at the boundary is unchanged, see Figure \ref{fig:billiard}. The radius of the circle is determined by the speed of the particle and the strength of the magnetic field. We fix the speed of the particle such that the radius is precisely the inverse of the strength of the magnetic field. We call this the Larmor radius.

Both billiards can be described as a map on $\Sigma\times(-\tfrac{\pi}{2},\tfrac{\pi}{2})$. The pair describes a point of incidence with outgoing direction. The billiard map sends such a pair to the next point of incidence together with the outgoing direction arising by following a straight line / fixed-radius circle. Thus we think of the billiard map as discrete dynamical system on $\Sigma\times(-\tfrac{\pi}{2},\tfrac{\pi}{2})$. It preserves the symplectic form $\om=\cos\phi \, d\phi\wedge ds $ where $s$ is the arc-length coordinate on $\Sigma$. Moreover, it is well-known that non-magnetic billiards with $\Sigma$ being an ellipse forms an integrable system, see for instance \cite{Tabachnikov_Geometry_and_billiards} and Example \ref{Fig:Ex0}. Recently, it was (roughly speaking) proved in \cite{Bialy_Mironov_Algebraic_non_integrability_of_magnetic_billiards} that the only (algebraic) integrable magnetic billiard occurs for $\Sigma$ being a circle, see the article for the precise statement.

KAM theory asserts that perturbations of the integrable billiard contains many invariant curves. We will numerically demonstrate this for perturbations being a non-zero magnetic field and a non-elliptical table.

\section{Numerical setup}

We briefly describe the geometric setup and the numerical algorithm. Our table $T$ is always bounded by the curve
\begin{align*}
\Sigma=\Big\{(x,y)\in\mathbb{R}^2\;:\; |x|^{p}+\frac{1}{1-\varepsilon^2}\,|y|^p=10^p\Big\}
\end{align*}
with  eccentricity $\varepsilon$ and power $p$ as free parameters, where $p=2$ corresponds to an ellipse. We simulate magnetic billiard inside $T$ choosing an external magnetic field of strength $B$ and assuming that the particle speed is normalized such that the Larmor radius is $B^{-1}$.  For parameterizing phase space we use angular coordinates $(\theta_\mathrm{pos},\theta_\mathrm{vel})\in(0,2\pi)\times(-\tfrac{\pi}{2},+\tfrac{\pi}{2})$ , where $\theta_{\mathrm{pos}}$ denotes the polar angle of points in $\Sigma$  and is hence not identical with the arc-length.

%%The numerical algorithm is illustrated in Figure \ref{Fig:FlowChart} and easy to implement.
%%For some choices of the parameters we computed a large number (= {\tt number of orbits}) of orbits as follows.  Initial data is chosen randomly, using the uniform probability distribution with respect to the variables $(\theta_\mathrm{pos},\theta_\mathrm{vel})\in(0,2\pi)\times(-\tfrac{\pi}{2},+\tfrac{\pi}{2})$, and each orbit undergoes a large number of reflections (= {\tt points per orbit}). 

\begin{figure}[ht!] %
\centering{ %
\includegraphics[width=0.9\textwidth]{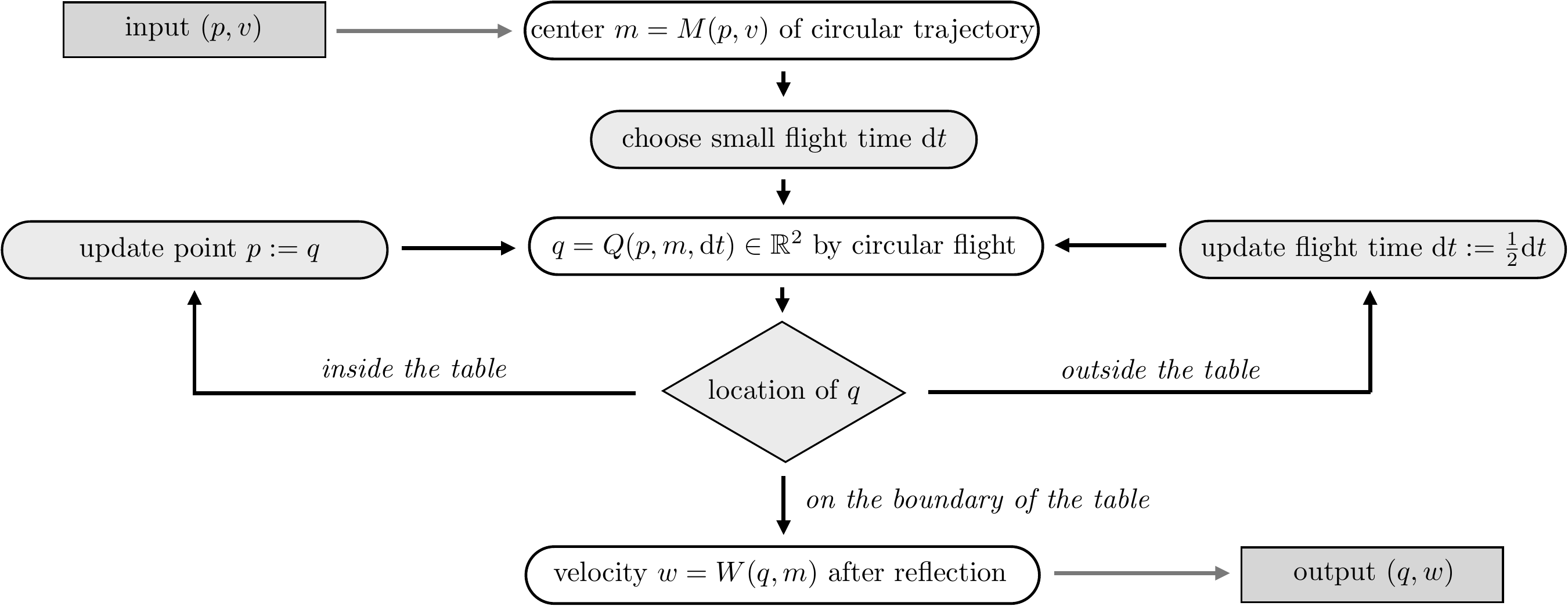}
}
\caption{Flow chart for the numerical computation of the billiard map with input and output belonging to $\Sigma\times(-\tfrac{\pi}{2},\tfrac{\pi}{2})$. The prescribed numerical accuracy of order $10^{-9}$ enters in determining whether the point $q$ lies on $\Sigma$.}
\label{Fig:FlowChart}
\end{figure}

The numerical algorithm is illustrated in Figure \ref{Fig:FlowChart} and easy to implement.
For some choices of the parameters we chose a large number (= {\tt number of orbits}) of random initial data and computed for each orbit a large number of particle reflections (= {\tt points per orbit}).  The underlying probability distribution was uniform with respect to  $(\theta_\mathrm{pos},\theta_\mathrm{vel})\in(0,2\pi)\times(-\tfrac{\pi}{2}+\delta,+\tfrac{\pi}{2}-\delta)$, where the small cut-off parameter $\delta$ excludes degenerate orbits.

The numerical results are displayed in the following figures. Each figure contains a phase space plot in which six orbits are colored. For these we also plot a certain number (= {\tt points on table}) of reflections in configuration space, i.e., on the table. 
\renewcommand{\figurename}{Example}
\setcounter{figure}{-1}
\begin{figure}
\centering{ %
\includegraphics[height=0.998\textwidth, angle=90]{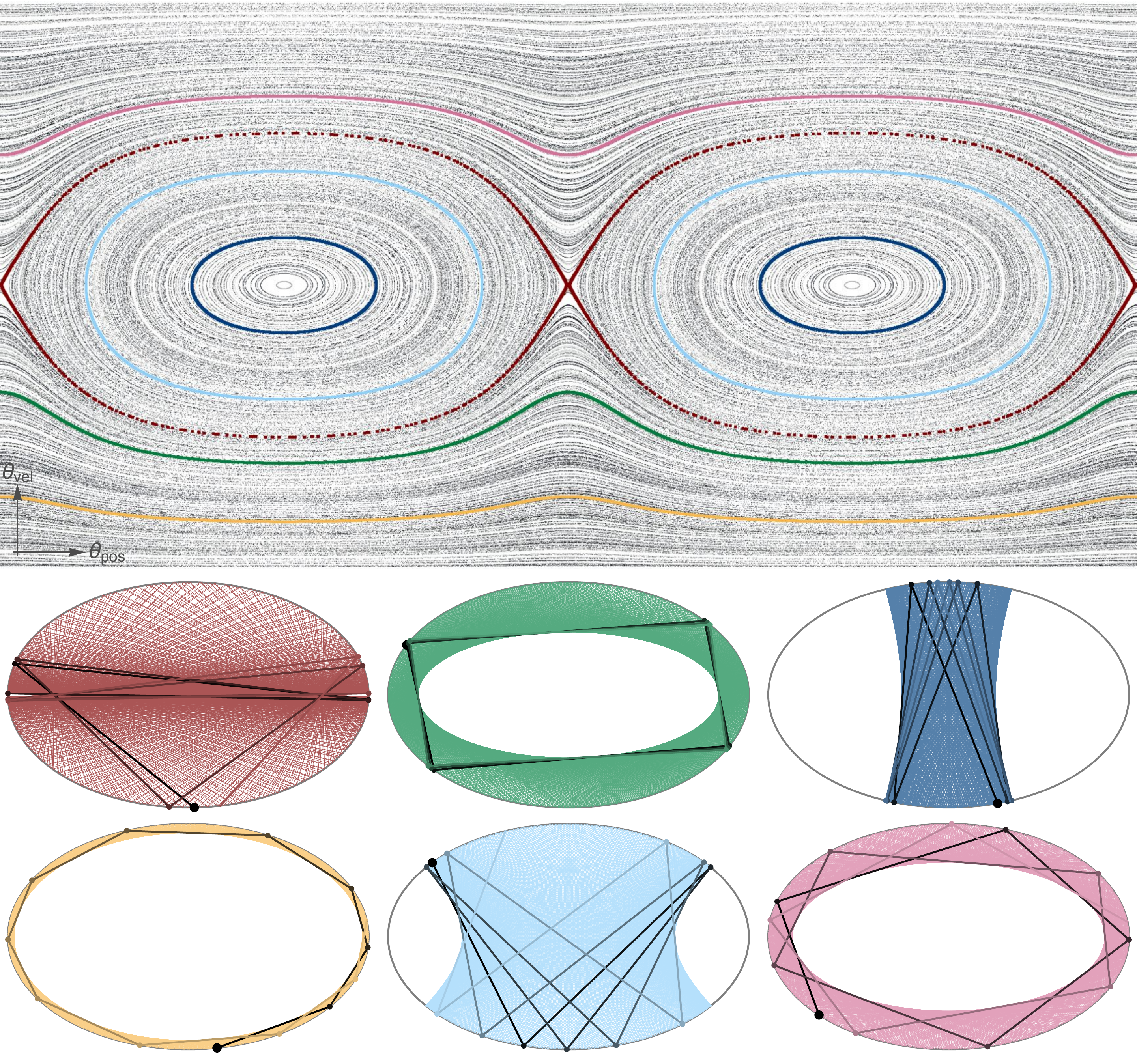}\\
\bigskip\bigskip\bigskip%
\boxed{ %
\begin{tabular}{lclclcl}
{\tt magnetic field} $B$ &=&0.0&&{\tt number of orbits}&=&1000\\
{\tt eccentricity} $\varepsilon$&=&1.5&&{\tt points per orbit}&=&1000\\
{\tt power} $p$ &=&2.0&&{\tt points on table}&=&1000\\
\end{tabular}
} }%
\bigskip\bigskip\bigskip%
\caption{Zero magnetic field and elliptic table.
}
\label{Fig:Ex0}
\end{figure}
\begin{figure}
\centering{ %
\includegraphics[height=0.998\textwidth, angle=90]{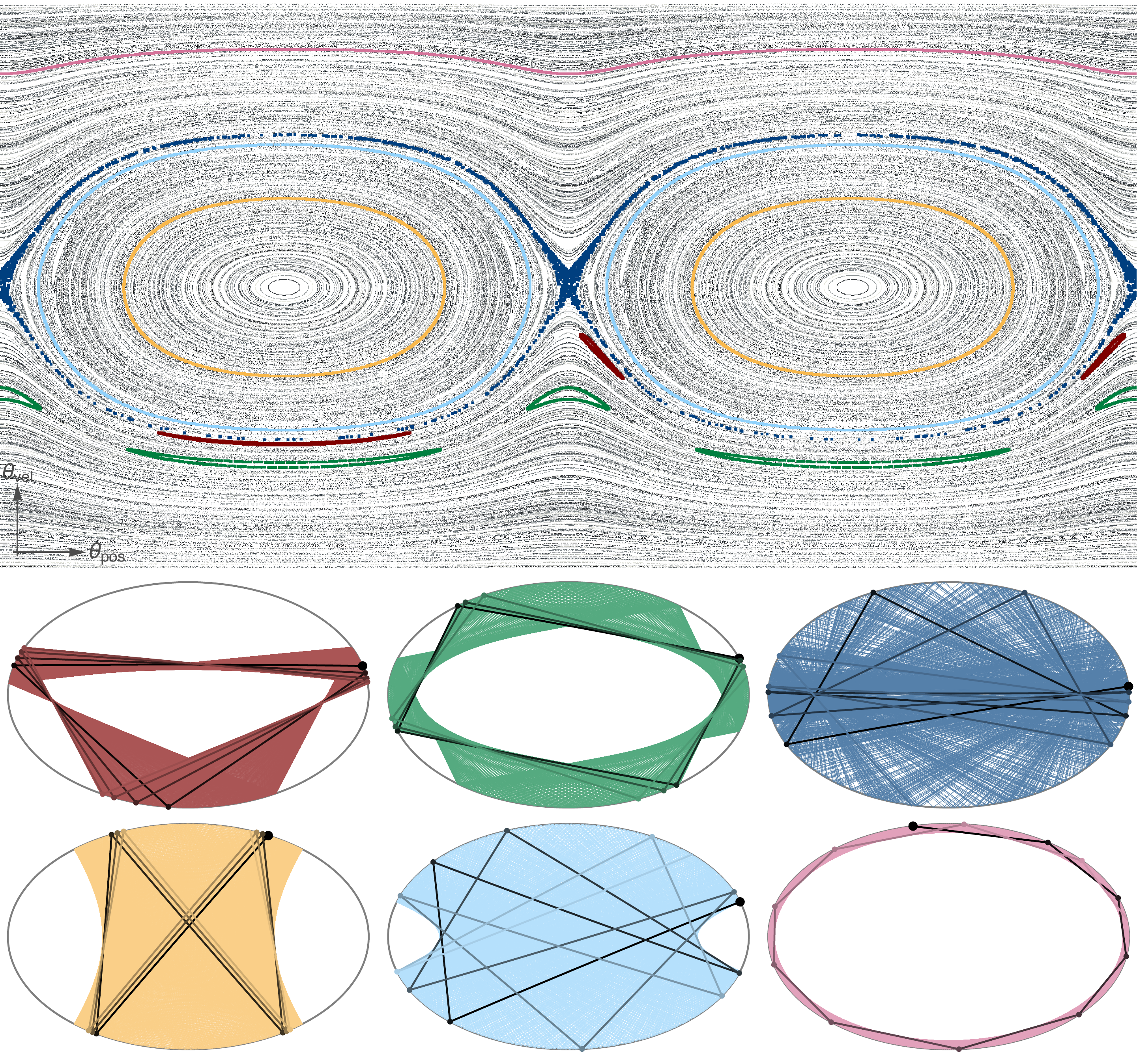}\\
\bigskip\bigskip\bigskip%
\boxed{ %
\begin{tabular}{lclclcl}
{\tt magnetic field} $B$ &=&0.01&&{\tt number of orbits}&=&1000\\
{\tt eccentricity} $\varepsilon$&=&1.5&&{\tt points per orbit}&=&1000\\
{\tt power} $p$ &=&2.0&&{\tt points on table}&=&2000\\
\end{tabular}
} }%
\bigskip\bigskip\bigskip%
\caption{Small magnetic field and elliptic table.
}
\label{Fig:Ex1}
\end{figure}
\begin{figure}
\centering{ %
\includegraphics[height=0.998\textwidth, angle=90]{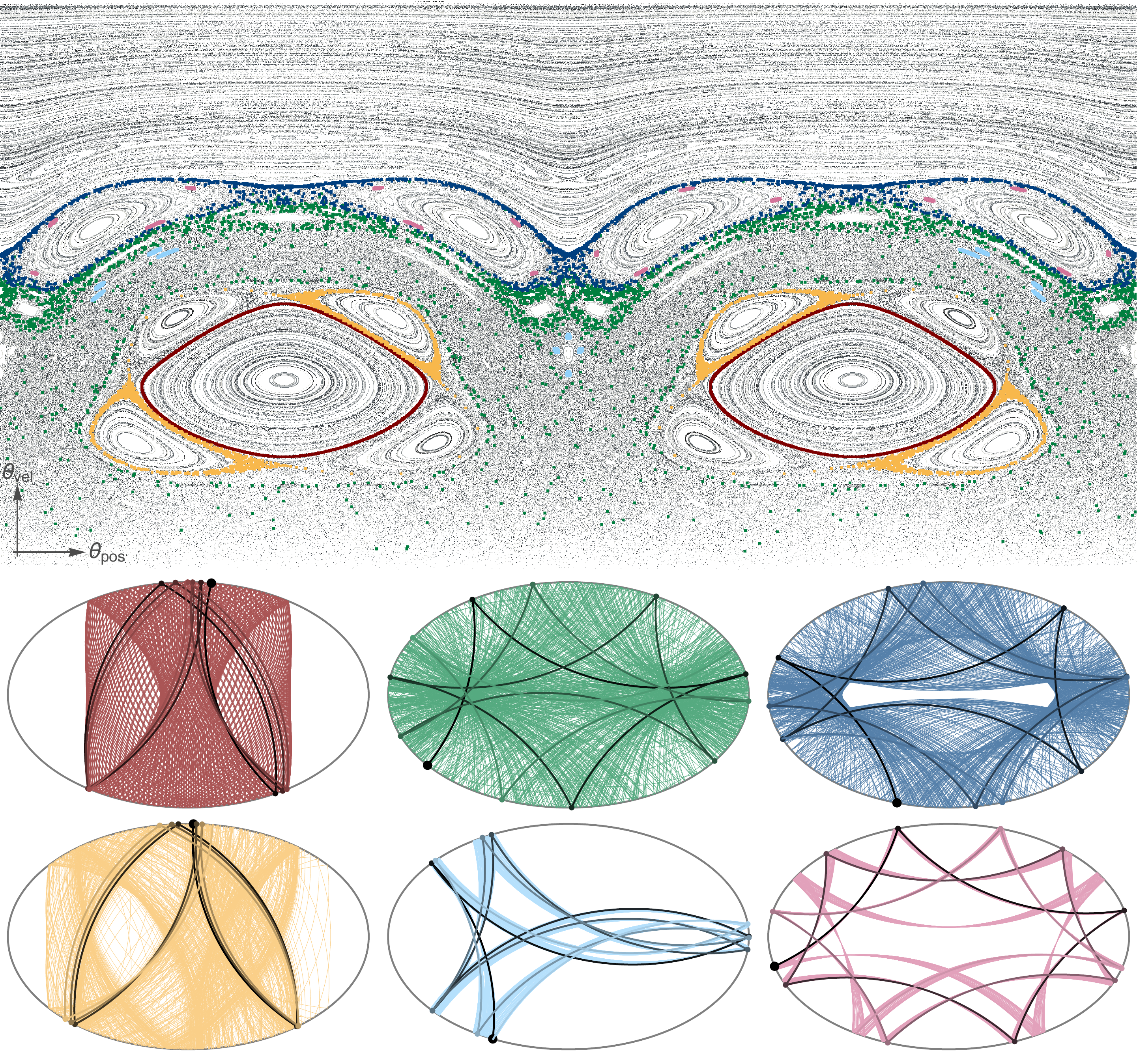}\\
\bigskip\bigskip\bigskip%
\boxed{ %
\begin{tabular}{lclclcl}
{\tt magnetic field} $B$&=&0.5&&{\tt number of orbits}&=&1000\\
{\tt eccentricity} $\varepsilon$&=&1.5&&{\tt points per orbit}&=&1000\\
{\tt power} $p$  &=&2.0&&{\tt points on table}&=&500\\
\end{tabular}
} }%
\bigskip\bigskip\bigskip%
\caption{Moderate magnetic field and elliptic table.
}
\label{Fig:Ex2}
\end{figure}

\begin{figure}
\centering{ %
\includegraphics[height=0.998\textwidth, angle=90]{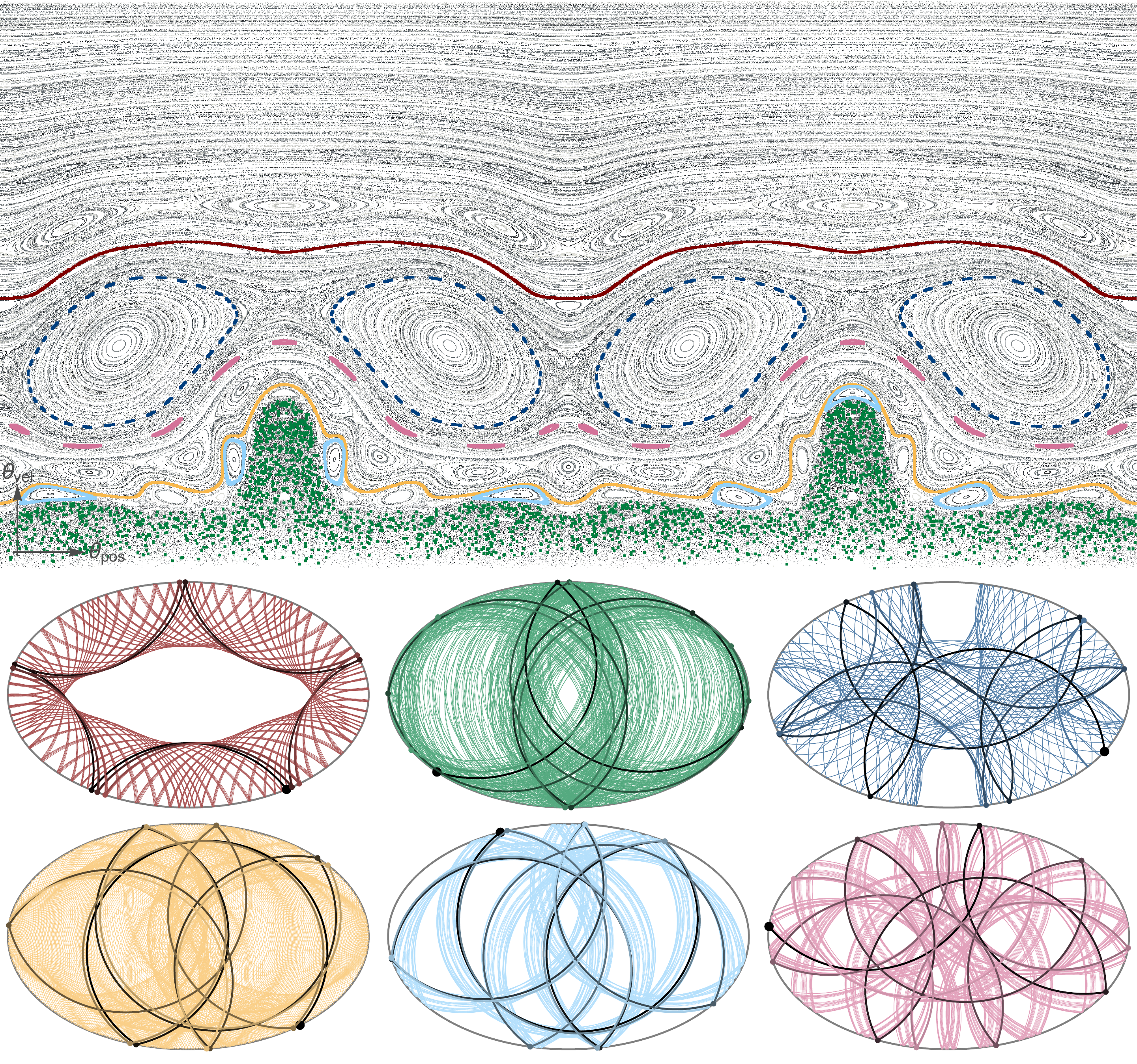}\\
\bigskip\bigskip\bigskip%
\boxed{ %
\begin{tabular}{lclclcl}
{\tt magnetic field} $B$&=&1.0&&{\tt number of orbits}&=&1000\\
{\tt eccentricity} $\varepsilon$ &=&1.5&&{\tt points per orbit}&=&1000\\
{\tt power} $p$ &=&2.0&&{\tt points on table}&=&500\\
\end{tabular}
} }%
\bigskip\bigskip\bigskip%
\caption{Large magnetic field and elliptic table.
}
\label{Fig:Ex3}
\end{figure}

\begin{figure}
\centering{ %
\includegraphics[height=0.998\textwidth, angle=90]{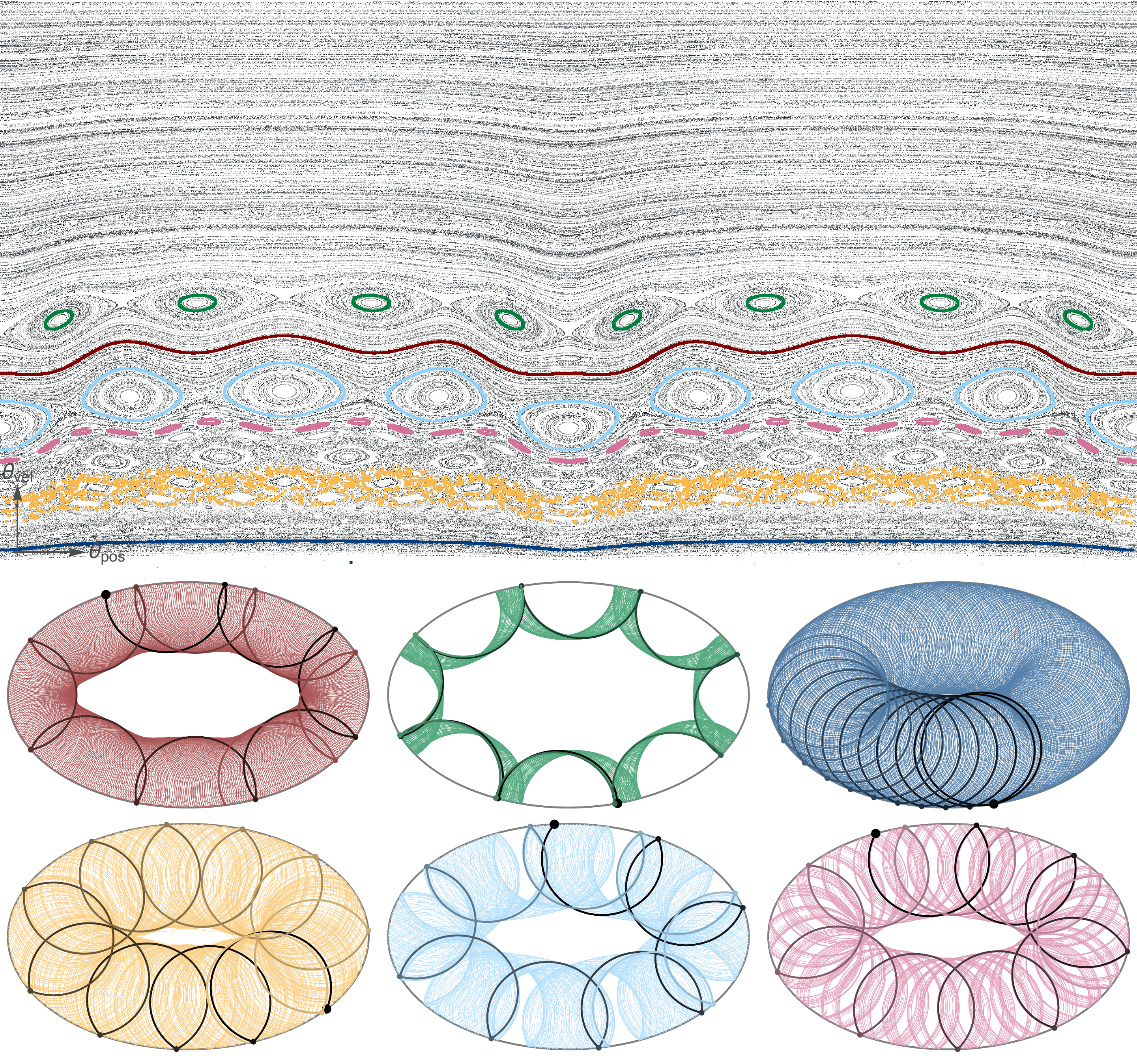}\\
\bigskip\bigskip\bigskip%
\boxed{\begin{tabular}{lclclcl}
{\tt magnetic field} $B$ &=&2.0&&{\tt number of orbits}&=&1000\\
{\tt eccentricity} $\varepsilon$ &=&1.5&&{\tt points per orbit}&=&1000\\
{\tt power} $p$ &=&2.0&&{\tt points on table}&=&500\\
\end{tabular}
} } %
\bigskip\bigskip\bigskip%
\caption{Even larger magnetic field and elliptic table.
}
\label{Fig:Ex4}
\end{figure}

\begin{figure}
\centering{ %
\includegraphics[height=0.998\textwidth, angle=90]{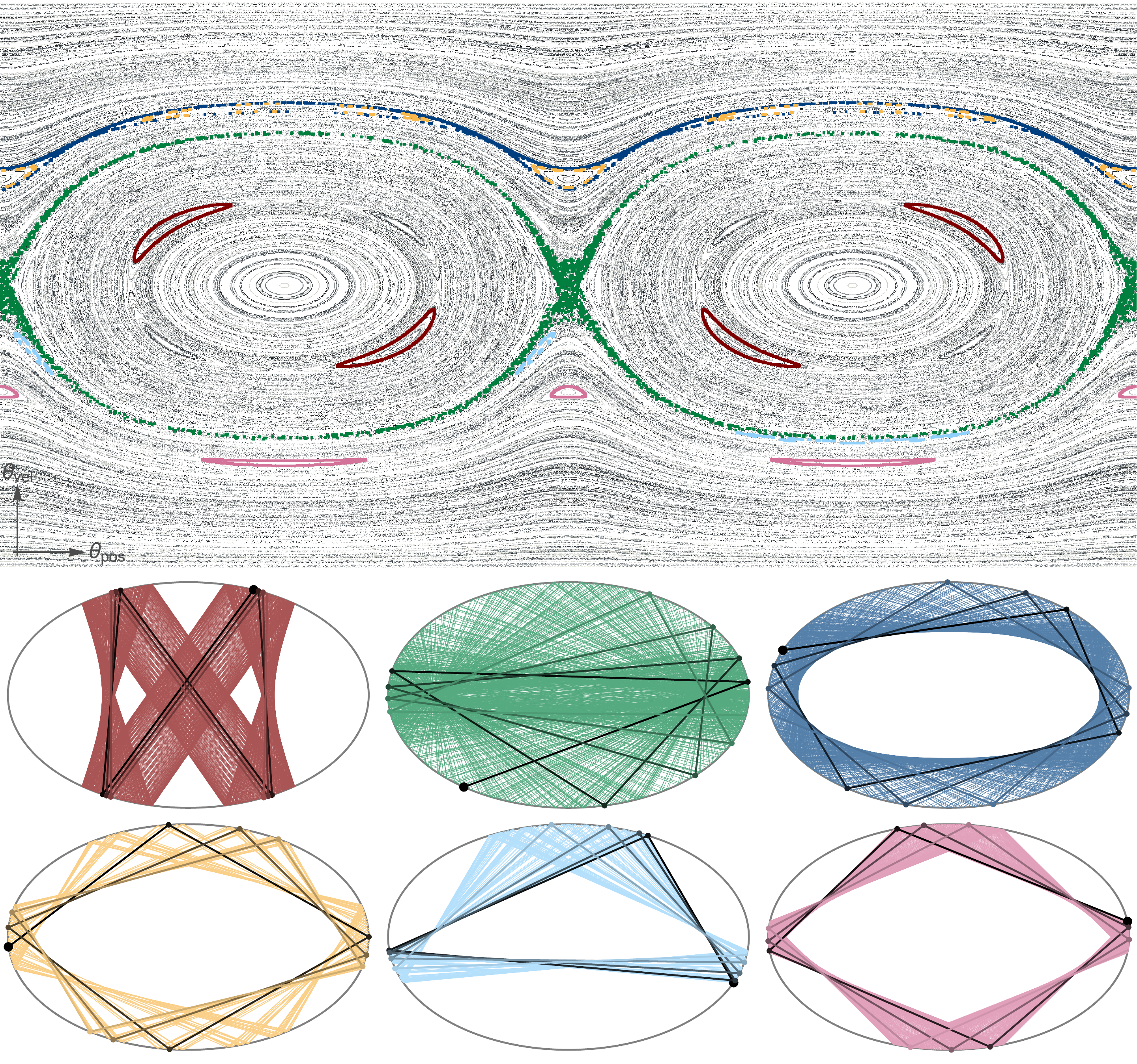}\\
\bigskip\bigskip\bigskip%
\boxed{\begin{tabular}{lclclcl}
{\tt magnetic field} $B$ &=&0.0&&{\tt number of orbits}&=&1000\\
{\tt eccentricity} $\varepsilon$ &=&1.5&&{\tt points per orbit}&=&1000\\
{\tt power} $p$ &=&2.005&&{\tt points on table}&=&1000\\
\end{tabular}
} }%
\bigskip\bigskip\bigskip%
\caption{Zero magnetic field and slightly non-elliptic table.
}
\label{Fig:Ex5}
\end{figure}

\begin{figure}
\centering{ %
\includegraphics[height=0.625\textwidth, angle=90]{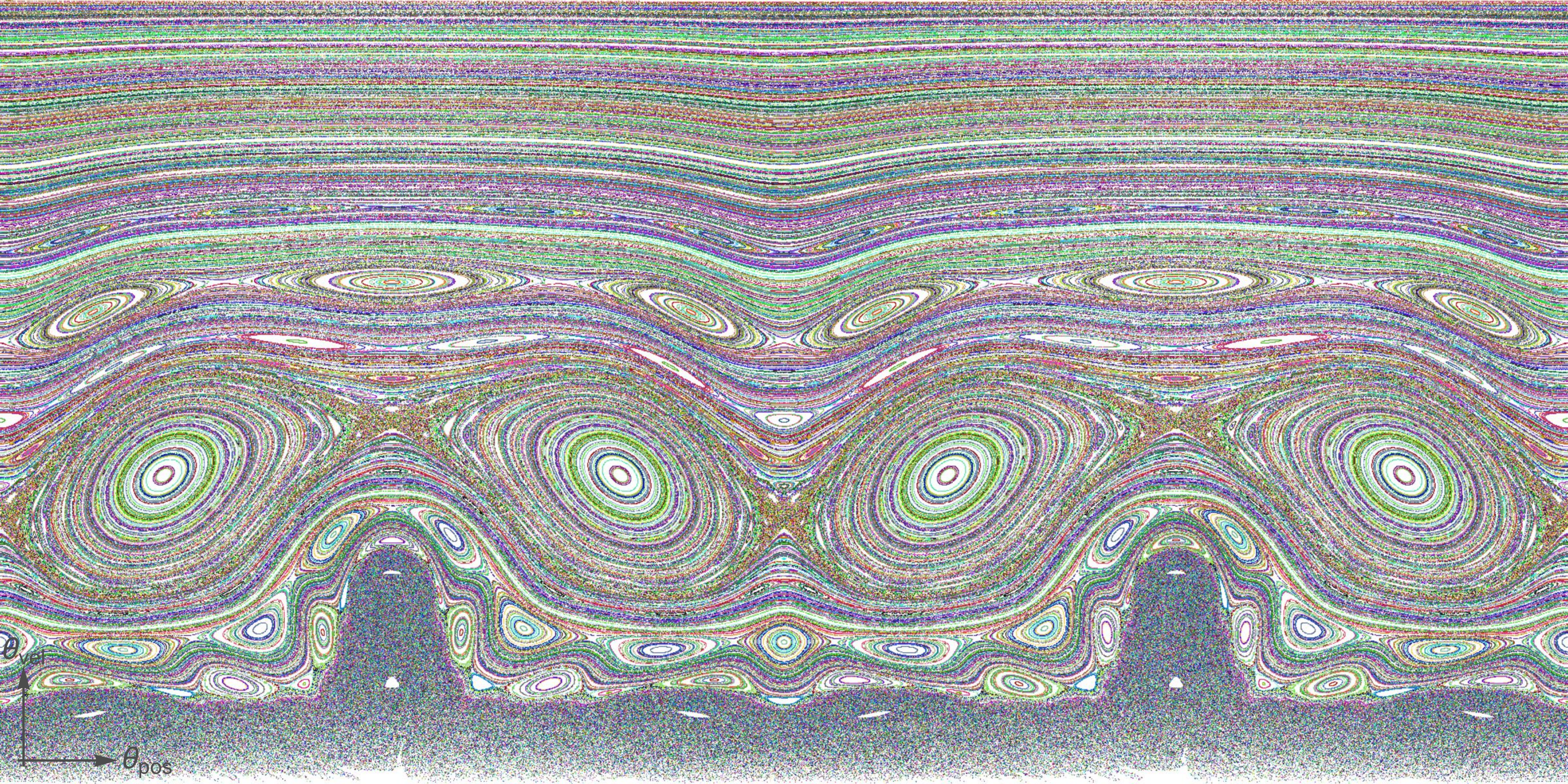}\\
\bigskip\bigskip\bigskip%
\boxed{ %
\begin{tabular}{lclclcl}
{\tt magnetic field} $B$ &=&1.0&&{\tt number of orbits}&=&2000\\
{\tt eccentricity} $\varepsilon$&=&1.5&&{\tt points per orbit}&=&3000\\
{\tt power} $p$ &=&2.0&&&&\\
\end{tabular}
} }%
\bigskip\bigskip\bigskip%
\caption{Phase portrait of Example 3 with higher resolution and random color for each orbit.
}
\label{Fig:All}
\end{figure}

%%%%%%%%%%%%%%%%%%%%%%%%%%%%%%%%%%%%%%%%%%%%%%%%%%%%%%%%%%%%%%%%%%%%%%%%%%%%%%%%%%%%%%%%%%%%%%%%%%%%%%%%%%%%%%%%%%%%%%%%%%%%%%%%%
\clearpage
\bibliographystyle{amsalpha}
\bibliography{bibtex_paper_list}

\end{document}